\catcode`\@=11

\font\tenmsx=msxm10
\font\sevenmsx=msxm7
\font\fivemsx=msxm5
\font\tenmsy=msym10
\font\sevenmsy=msym7
\font\fivemsy=msym5
\newfam\msxfam
\newfam\msyfam
\textfont\msxfam=\tenmsx  \scriptfont\msxfam=\sevenmsx
  \scriptscriptfont\msxfam=\fivemsx
\textfont\msyfam=\tenmsy  \scriptfont\msyfam=\sevenmsy
  \scriptscriptfont\msyfam=\fivemsy

\def\hexnumber@#1{\ifnum#1<10 \number#1\else
 \ifnum#1=10 A\else\ifnum#1=11 B\else\ifnum#1=12 C\else
 \ifnum#1=13 D\else\ifnum#1=14 E\else\ifnum#1=15 F\fi\fi\fi\fi\fi\fi\fi}

\def\msx@{\hexnumber@\msxfam}
\def\msy@{\hexnumber@\msyfam}
\mathchardef\boxdot="2\msx@00
\mathchardef\boxplus="2\msx@01
\mathchardef\boxtimes="2\msx@02
\mathchardef\square="0\msx@03
\mathchardef\blacksquare="0\msx@04
\mathchardef\centerdot="2\msx@05
\mathchardef\lozenge="0\msx@06
\mathchardef\blacklozenge="0\msx@07
\mathchardef\circlearrowright="3\msx@08
\mathchardef\circlearrowleft="3\msx@09
\mathchardef\rightleftharpoons="3\msx@0A
\mathchardef\leftrightharpoons="3\msx@0B
\mathchardef\boxminus="2\msx@0C
\mathchardef\Vdash="3\msx@0D
\mathchardef\Vvdash="3\msx@0E
\mathchardef\vDash="3\msx@0F
\mathchardef\twoheadrightarrow="3\msx@10
\mathchardef\twoheadleftarrow="3\msx@11
\mathchardef\leftleftarrows="3\msx@12
\mathchardef\rightrightarrows="3\msx@13
\mathchardef\upuparrows="3\msx@14
\mathchardef\downdownarrows="3\msx@15
\mathchardef\upharpoonright="3\msx@16

\mathchardef\downharpoonright="3\msx@17
\mathchardef\upharpoonleft="3\msx@18
\mathchardef\downharpoonleft="3\msx@19
\mathchardef\rightarrowtail="3\msx@1A
\mathchardef\leftarrowtail="3\msx@1B
\mathchardef\leftrightarrows="3\msx@1C
\mathchardef\rightleftarrows="3\msx@1D
\mathchardef\Lsh="3\msx@1E
\mathchardef\Rsh="3\msx@1F
\mathchardef\rightsquigarrow="3\msx@20
\mathchardef\leftrightsquigarrow="3\msx@21
\mathchardef\looparrowleft="3\msx@22
\mathchardef\looparrowright="3\msx@23
\mathchardef\circeq="3\msx@24
\mathchardef\succsim="3\msx@25
\mathchardef\gtrsim="3\msx@26
\mathchardef\gtrapprox="3\msx@27
\mathchardef\multimap="3\msx@28
\mathchardef\therefore="3\msx@29
\mathchardef\because="3\msx@2A
\mathchardef\doteqdot="3\msx@2B

\mathchardef\triangleq="3\msx@2C
\mathchardef\precsim="3\msx@2D
\mathchardef\lesssim="3\msx@2E
\mathchardef\lessapprox="3\msx@2F
\mathchardef\eqslantless="3\msx@30
\mathchardef\eqslantgtr="3\msx@31
\mathchardef\curlyeqprec="3\msx@32
\mathchardef\curlyeqsucc="3\msx@33
\mathchardef\preccurlyeq="3\msx@34
\mathchardef\leqq="3\msx@35
\mathchardef\leqslant="3\msx@36
\mathchardef\lessgtr="3\msx@37
\mathchardef\backprime="0\msx@38
\mathchardef\risingdotseq="3\msx@3A
\mathchardef\fallingdotseq="3\msx@3B
\mathchardef\succcurlyeq="3\msx@3C
\mathchardef\geqq="3\msx@3D
\mathchardef\geqslant="3\msx@3E
\mathchardef\gtrless="3\msx@3F
\mathchardef\sqsubset="3\msx@40
\mathchardef\sqsupset="3\msx@41
\mathchardef\trianglerighteq="3\msx@44
\mathchardef\trianglelefteq="3\msx@45
\mathchardef\bigstar="0\msx@46
\mathchardef\between="3\msx@47
\mathchardef\blacktriangledown="0\msx@48
\mathchardef\blacktriangleright="3\msx@49
\mathchardef\blacktriangleleft="3\msx@4A
\mathchardef\blacktriangle="0\msx@4E
\mathchardef\triangledown="0\msx@4F
\mathchardef\eqcirc="3\msx@50
\mathchardef\lesseqgtr="3\msx@51
\mathchardef\gtreqless="3\msx@52
\mathchardef\lesseqqgtr="3\msx@53
\mathchardef\gtreqqless="3\msx@54
\mathchardef\Rrightarrow="3\msx@56
\mathchardef\Lleftarrow="3\msx@57
\mathchardef\veebar="2\msx@59
\mathchardef\barwedge="2\msx@5A
\mathchardef\doublebarwedge="2\msx@5B
\mathchardef\angle="0\msx@5C
\mathchardef\measuredangle="0\msx@5D
\mathchardef\sphericalangle="0\msx@5E
\mathchardef\varpropto="3\msx@5F
\mathchardef\smallsmile="3\msx@60
\mathchardef\smallfrown="3\msx@61
\mathchardef\Subset="3\msx@62
\mathchardef\Supset="3\msx@63
\mathchardef\Cup="2\msx@64

\mathchardef\Cap="2\msx@65

\mathchardef\curlywedge="2\msx@66
\mathchardef\curlyvee="2\msx@67
\mathchardef\leftthreetimes="2\msx@68
\mathchardef\rightthreetimes="2\msx@69
\mathchardef\subseteqq="3\msx@6A
\mathchardef\supseteqq="3\msx@6B
\mathchardef\bumpeq="3\msx@6C
\mathchardef\Bumpeq="3\msx@6D
\mathchardef\lll="3\msx@6E

\mathchardef\ggg="3\msx@6F

\mathchardef\circledS="0\msx@73
\mathchardef\pitchfork="3\msx@74
\mathchardef\dotplus="2\msx@75
\mathchardef\backsim="3\msx@76
\mathchardef\backsimeq="3\msx@77
\mathchardef\complement="0\msx@7B
\mathchardef\intercal="2\msx@7C
\mathchardef\circledcirc="2\msx@7D
\mathchardef\circledast="2\msx@7E
\mathchardef\circleddash="2\msx@7F
\def\ulcorner{\delimiter"4\msx@70\msx@70 }
\def\urcorner{\delimiter"5\msx@71\msx@71 }
\def\llcorner{\delimiter"4\msx@78\msx@78 }
\def\lrcorner{\delimiter"5\msx@79\msx@79 }
\def\yen{\mathhexbox\msx@55 }
\def\checkmark{\mathhexbox\msx@58 }
\def\circledR{\mathhexbox\msx@72 }
\def\maltese{\mathhexbox\msx@7A }
\mathchardef\lvertneqq="3\msy@00
\mathchardef\gvertneqq="3\msy@01
\mathchardef\nleq="3\msy@02
\mathchardef\ngeq="3\msy@03
\mathchardef\nless="3\msy@04
\mathchardef\ngtr="3\msy@05
\mathchardef\nprec="3\msy@06
\mathchardef\nsucc="3\msy@07
\mathchardef\lneqq="3\msy@08
\mathchardef\gneqq="3\msy@09
\mathchardef\nleqslant="3\msy@0A
\mathchardef\ngeqslant="3\msy@0B
\mathchardef\lneq="3\msy@0C
\mathchardef\gneq="3\msy@0D
\mathchardef\npreceq="3\msy@0E
\mathchardef\nsucceq="3\msy@0F
\mathchardef\precnsim="3\msy@10
\mathchardef\succnsim="3\msy@11
\mathchardef\lnsim="3\msy@12
\mathchardef\gnsim="3\msy@13
\mathchardef\nleqq="3\msy@14
\mathchardef\ngeqq="3\msy@15
\mathchardef\precneqq="3\msy@16
\mathchardef\succneqq="3\msy@17
\mathchardef\precnapprox="3\msy@18
\mathchardef\succnapprox="3\msy@19
\mathchardef\lnapprox="3\msy@1A
\mathchardef\gnapprox="3\msy@1B
\mathchardef\nsim="3\msy@1C
\mathchardef\napprox="3\msy@1D
\mathchardef\nsubseteqq="3\msy@22
\mathchardef\nsupseteqq="3\msy@23
\mathchardef\subsetneqq="3\msy@24
\mathchardef\supsetneqq="3\msy@25
\mathchardef\subsetneq="3\msy@28
\mathchardef\supsetneq="3\msy@29
\mathchardef\nsubseteq="3\msy@2A
\mathchardef\nsupseteq="3\msy@2B
\mathchardef\nparallel="3\msy@2C
\mathchardef\nmid="3\msy@2D
\mathchardef\nshortmid="3\msy@2E
\mathchardef\nshortparallel="3\msy@2F
\mathchardef\nvdash="3\msy@30
\mathchardef\nVdash="3\msy@31
\mathchardef\nvDash="3\msy@32
\mathchardef\nVDash="3\msy@33
\mathchardef\ntrianglerighteq="3\msy@34
\mathchardef\ntrianglelefteq="3\msy@35
\mathchardef\ntriangleleft="3\msy@36
\mathchardef\ntriangleright="3\msy@37
\mathchardef\nleftarrow="3\msy@38
\mathchardef\nrightarrow="3\msy@39
\mathchardef\nLeftarrow="3\msy@3A
\mathchardef\nRightarrow="3\msy@3B
\mathchardef\nLeftrightarrow="3\msy@3C
\mathchardef\nleftrightarrow="3\msy@3D
\mathchardef\divideontimes="2\msy@3E
\mathchardef\varnothing="0\msy@3F
\mathchardef\nexists="0\msy@40
\mathchardef\mho="0\msy@66
\mathchardef\thorn="0\msy@67
\mathchardef\beth="0\msy@69
\mathchardef\gimel="0\msy@6A
\mathchardef\daleth="0\msy@6B
\mathchardef\lessdot="3\msy@6C
\mathchardef\gtrdot="3\msy@6D
\mathchardef\ltimes="2\msy@6E
\mathchardef\rtimes="2\msy@6F
\mathchardef\shortmid="3\msy@70
\mathchardef\shortparallel="3\msy@71
\mathchardef\smallsetminus="2\msy@72
\mathchardef\thicksim="3\msy@73
\mathchardef\thickapprox="3\msy@74
\mathchardef\approxeq="3\msy@75
\mathchardef\succapprox="3\msy@76
\mathchardef\precapprox="3\msy@77
\mathchardef\curvearrowleft="3\msy@78
\mathchardef\curvearrowright="3\msy@79
\mathchardef\digamma="0\msy@7A
\mathchardef\varkappa="0\msy@7B
\mathchardef\hslash="0\msy@7D
\mathchardef\hbar="0\msy@7E
\mathchardef\backepsilon="3\msy@7F
\def\Bbb{\ifmmode\let\next\Bbb@\else
 \def\next{\errmessage{Use \string\Bbb\space only in math mode}}\fi\next}
\def\Bbb@#1{{\Bbb@@{#1}}}
\def\Bbb@@#1{\fam\msyfam#1}

\catcode`\@=12
\def\C{\Bbb C} 
\def\R{\Bbb R} 
\def\P{\Bbb P} 
\def\Z{\Bbb Z}

\documentstyle{article}

\begin{document}
\sloppy
\newtheorem{Th}{Theorem}[section]
\newtheorem{Satz}{Satz}[section]
\newtheorem{Prop}{Proposition}[section]
\newtheorem{Lemma}{Lemma}[section]
\newtheorem{Rem}{Remark}[section]
\newtheorem{Def}{Definition}[section]
\newtheorem{Cor}{Corollary}[section]
\title{On the automorphism groups of algebraic bounded domains.}
\author{D.~Zaitsev}
\maketitle


\tableofcontents
\parskip 1.5mm
	\section{Introduction}\label{in}

	Let  $D$  be  a  domain in $\C^n$ and $Aut(D)$ be the group of
all biholomorphic automorphisms of $D$. Let $v\in D$ be fixed and define
the map $C_v\colon Aut(D) \to D\times Gl(n)$ by $f\mapsto (f(v),f_{*v})$.
The theorem of H.~Cartan (see Narasimhan, \cite{N}, p.~169) can be stated
as follows.

\begin{samepage}
\begin{Th}
Let $D$ be bounded. Then:
\begin{enumerate}
\item The group $Aut(D)$ possesses a natural Lie group structure
compatible with the compact-open topology such that the action
$Aut(D)\times D\to D$ is real analytic.
\item For all $v\in D$ the map $C_v$ is a real-analytic homeomorphism
onto its image.
\end{enumerate}
\end{Th}
\end{samepage}

	The domains we discuss here are open connected sets defined by finitely
many real polynomial inequalities or connected finite unions of such sets.
These are the domains in the so-called ``semi-algebraic category''
defined below (Definition~\ref{def-s-a}).
For this reason we call them ``semi-algebraic domains''.
In this paper we are interested in the algebraic nature of the image of
$Aut(D)$ and its subgroups in $D\times Gl(n)$ under the map $C_v$.

	{\bf Example.1.} The simplest example of a semi-algebraic domain is the
unit disk $D=\{|z|<1\}$. For this domain
$$Aut(D)=PGL_2(\R)_+:= \{A\in PGL_2(\R) \mid \det A > 0 \}.$$
We see that $Aut(D)$, as a subgroup of $PGL_2(\R)$, is defined by an
inequality and therefore is not an algebraic subgroup. In fact, the
group $Aut(D)$ here does not admit algebraic structure as a Lie group.

	To show this, assume that there is a Lie isomorphism
$\varphi\colon PGL_2(\R)_+\to G$, where $G$ is a real algebraic group.
It continues to an isomorphism of complexifications
$\varphi^{\C}\colon PGL_2(\C) \to G^{\C}$. The latter, being a Lie isomorphism
between semi-simple complex algebraic groups, is algebraic. Since
$G\subset G^{\C}$ is real algebraic, so is its preimage
$(\varphi^{\C})^{-1}(G)=PGL_2(\R)_+$. On the other hand,
$PGL_2(\R)_+\subset PGL_2(\C)$ is not real Zariski closed. This is a
contradiction.
\par\hfill {\bf Q. E. D.}

	{\bf Example.2.} More generally let $D$ be a bounded homogeneous
domain in $\C^n$. By the classification theorem of Vinberg, Gindikin
and Pyatetskii-Shapiro (see \cite{VGP}, Theorem~6, p.~434), $D$ is
biholomorphic to a homogeneous Siegel domain of the 1st or the 2nd kind.
Such a domain is defined algebraically in terms
of a homogeneous convex cone (\cite{VGP}) .
Rothaus (\cite{R}) gave a procedure for constructing
all homogeneous convex cones. The construction implies that
all homogeneous convex cones, and therefore all homogeneous Siegel domains
of the 1st and 2nd kind, are defined by finitely many polynomial inequalities.
The Siegel domains are unbounded but they are birationally equivalent to
bounded domains which are also defined by finitely many polynomial inequalities
and therefore are semi-algebraic. Thus, $D$ is biholomorphic to a bounded
semi-algebaic domain.

	The automorphism group $Aut(D)$ of a bounded homogeneous domain was
discussed by Kaneyuki (see \cite{K}) where he proved in particular
that the identity component $Aut(D)^0$ is isomorphic to an identity
component of a real algebraic group (\cite{K}, Theorem.3.2., p.106).
Let $x_0\in D$ be a fixed point. Since $D=Aut(D)/Iso(x_0)$
and the isotropy group $Iso(x_0)$
is compact, the automorphism group $Aut(D)$ has finitely many components.
Together with the result of Kaneyuki this implies that $Aut(D)$ is
isomorphic to an open subgroup of a real algebraic group.

	In the above examples the automorphism group $Aut(D)$ is isomorphic to
an open subgroup of a real algebraic group. Therefore, it admits a faithful
representation. For general domains however the automorphism group does not
admit a faithfull representation.

	{\bf Example.3.} Let
$$D:=\{(z,w)\in\C^2 \mid |z|^2+|w|^2<1, w\ne 0 \} $$
be the unit ball in $\C^2$ with a unit disk removed. We consider $D$ as
a subset of $\P^2$ with homogeneous coordinates $\xi_0,\xi_1,\xi_2$,
$z=\xi_1/\xi_0$, $w=\xi_2/\xi_0$. The automorphism group of $D$
is $Aut(D)=SU(1,1)\times S^1$ with the action on $D$ given by
$$ (A,\tau)(\xi_0,\xi_1,\xi_2) = \pmatrix{
A       & 0 \atop 0 \cr
0 \, 0  & \tau      \cr} \pmatrix{
\xi_0 \cr
\xi_1 \cr
\xi_2 \cr },\quad A\in SU(1,1),\quad \tau \in S^1. $$
The isomorphism $K\colon C\mapsto JCJ^{-1}$ between $SL_2(\R)$ and
$SU(1,1)$ ($J=\pmatrix{-i&1\cr i&1\cr}$) yelds an effective action of
$SL_2(\R)$ on $D$. Furthermore, each map $j_x\colon SL_2(\R)\to D$,
$j_x(C):=Cx$ induces an isomorphism between fundamental groups
$\pi_1(SL_2(\R))=\pi_1(D)=\Z$.
This implies that the induced action of a finite covering of $SL_2(\R)$
on $D$ lifts to an effective action on the finite covering $D'$ of $D$ of the
same degree. But no covering of $SL_2(\R)$ admits a faithfull
representation (because every such representation factorizes through a
representation of $SL_2(\C)$). On the other hand, the map
$(z,w)\mapsto (z,\root d \of w)$ defines
an isomorphism between a finite covering of $D$ of degree $d$ and a bounded
semi-algebraic domain $\tilde D\subset\C^2$.

	In this example the group $Aut(D)$ is not isomorphic to an open subset
of an algebraic group. Moreover, even in case it is,
the action $Aut(D)\times D\to D$ can be ``far from algebraic''.
This phenomena is shown in the following example.

	{\bf Example.4.} Let $F=\C/\Lambda$ be a complex elliptic curve and
${\cal P}\colon F\to \P^1$ the Weierstra\ss\ ${\cal P}$-function which
defines a $2-1$ ramified covering over $P^1$. The strip
$\{z\in\C \mid c-\epsilon<{\rm Im}z<c+\epsilon\}$ covers a ``circle strip''
$\tilde D\subset \C/\Lambda$. Let $D$ be the projection of $\tilde D$
on $P^1$. If the constant $c$ is generic and $\epsilon$ is small enough,
the projection of $\tilde D$ is biholomorphic. The real algebraic group
$S^1\subset \C/\Lambda$ acts by translations on $\tilde D$ which yields
an effective action on $D$. Since the domain $D$ is bounded by real elliptic
curves,  it is therefore semi-algebraic. However the action of $S^1$ is
expressed in terms of the Weierstra\ss\ ${\cal P}$-function and is not
algebraic. In fact there is no homomorphism of $S^1$ in a real algebraic
group $G$ such that the action $S^1\times D\to D$ is given by restrictions
of polynomials on $G$. This follows from the classification of 1-dimensional
Nash (semi-algebraic) groups given by Madden and Stanton (see \cite{MS}).

	We see therefore that even in simple cases the class of
real algebraic groups and their subgroups is not large enough
to describe the group $Aut(D)$ and its action on $D$.
Consequently, we consider a larger class of groups where {\em defining
inequalities} are allowed.

\begin{Def}
\begin{enumerate}
\item A {\bf Nash function} is a real analytic function
$f=(f_1,\ldots,f_m)\colon U\to \R^m$ (where $U$ is an open semi-algebraic
subset of $\R^n$) such that for each of the components $f_k$ there is a
nontrivial polynomial $P$ with $P(x_1,\ldots,x_n,f_k(x_1,\ldots,x_n))=0$
for all $(x_1,\ldots,x_n)\in U$.
\item A {\bf Nash manifold} $M$ is a real analytic manifold with finitely many
coordinate charts $\phi_i\colon U_i\to V_i$ such that $V_i\subset\R^n$ is
semi-algebraic for all $i$ and the transition functions are Nash
(a Nash atlas).
\item A Nash manifold is called {\bf affine} if it can be Nash
(locally closed) imbedded in $\R^N$ for some $N$.
\item A {\bf Nash group} is a Nash manifold with a group operation
$(x,y)\to xy^{-1}$ which is Nash with respect to every Nash coordinate chart.
\end{enumerate}
\end{Def}

\begin{Rem}
	The simplest example of a Nash manifold which is not affine is
the quotient $\R/\Z$ with the Nash structure inherited from the standard
Nash structure on $\R$.
For the classification of such groups in the one-dimensional case see
J.~J.~Madden and C.~M.~Stanton in \cite{MS}.
\end{Rem}

	Roughly speaking, the goal of this paper is to prove that the automorphism
group $Aut(D)$ of a semi-algebraic domain $D$ has a natural Nash group
structure such that the action $Aut(D)\times D\to D$ is also Nash. For this
we need a certain non-degeneracy condition on the boundary of $D$.
To give the reader a flavour of the main result, we first mention
an application for the {\it algebraic} domains introduced by Diederich and
Forn\ae ss (\cite{DF}), for which this condition is automatically satisfied.

\begin{Def}(see Diederich-Forn\ae ss, \cite{DF})
	A domain $D\subset\subset C^n$ is called {\bf algebraic} if there exists
a real polynomial $r(z,\bar z)$ such that $D$ is a connected component of the
set
$$\{z\in\C^n \mid r(z,\bar z)<0 \}$$
and $dr(z)\ne 0$ for $z\in\partial D$.
\end{Def}

\begin{Th}\label{cor}
	Let $D\subset\subset\C^n$, $n>1$, be an algebraic domain.
The group $Aut(D)$ possesses a unique structure of an affine Nash group
so that the action $Aut(D)\times D\to D$ is Nash.
For all $v\in D$, $C_v\colon Aut(D)\to D\times Gl(n)$ is a Nash
isomorphism onto its image.
\end{Th}

\begin{Cor}
	Let $D$ be as in Theorem~\ref{cor}. Then the group $Aut(D)$
has finitely many connected components.
\end{Cor}

\begin{Rem}
	In general the number of components of $Aut(D)$ can be infinite.
For example, let $H:=\{z\in\C \mid {\rm Im}z=0\}$ be the upper half-plane and
$$\hat H:=H\cup \bigcup_{n\in\Z} B_{\epsilon}(n),$$
where $B_{\epsilon}(n)$ is the ball with centre $n$ and radius
$\epsilon<1/2$. Let $D\subset\C^2$ be the union of $H\times H$ and
$(H+i)\times \hat H$. Then $D$ is biholomorphic to a simply connected
bounded domain. The flat pieces of the boundary of $D$ admit canonical
foliations $z=\rm const$ and $w=\rm const$. The latters induce foliations
of $D$ of the same form which are preserved by the automorphisms
(see Remmert and Stein \cite{RS}).
By this argument one shows that $Aut(D)=\R\oplus\Z$.
\end{Rem}

For the formulation of our main result we need the following condition on
$D$, which is automatically satisfied for all bounded domains with smooth
boundary, in particular, for all algebraic domains.

\begin{Def}\label{L}
\begin{enumerate}
\item
	A domain $D\subset\C^n$ or its boundary
is called {\bf Levi-non-degenerate},
if there exists a point $x_0\in\partial D$ and a neighborhood $U\subset\C^n$
of $x_0$ such that
$$D\cap U=\{z\in U \mid \varphi(z)<0 \}$$
for a $C^2$-function
$\varphi$ with $d\varphi\ne 0$ and such that the Levi form in $x_0$
$$L_r(x_0):=\sum_{k,l=1}^n
{\partial^2 r \over \partial z_i \, \partial \bar z_j} dz_k \otimes d\bar z_l$$
restricted to the holomorphic tangent space of $\partial D$ in $x_0$ is
non-degenerate;
\item
	A domain $D\subset\C^n$ or its boundary
is called {\bf completely Levi-non-degenerate},
if every boundary point outside a real analytic subset of dimension $2n-2$
is non-degenerate in the above sense.
\end{enumerate}
\end{Def}

   Now let $D$ be semi-algebraic and
consider the subset $Aut_a(D)\subset Aut(D)$ of all
(biholomorphic) automorphisms which are Nash.
In the following sense
Nash automorphisms are branches of algebraic maps
(see Proposition~\ref{bra}).

\begin{Def}\label{br}
	Let $D$ be a domain in $\C^n$. A holomorphic map $f\in Aut(D)$
is a {\bf branch of an algebraic map} if
there exists a complex $n$-dimensional algebraic subvariety
$G\subset\C^n\times\C^n$ which contains the graph of $f$.
\end{Def}

   The following is the main result of this paper.

\begin{Th}\label{main}
	Let $D\subset\subset\C^n$, $n>1$, be a semi-algebraic
Levi-non-degenerate domain. Then
\begin{enumerate}
\item $Aut_a(D)$ is a (closed) Lie subgroup of $Aut(D)$,
\item $Aut_a(D)$ possesses a unique structure of an affine Nash group
so that the action $Aut_a(D)\times D\to D$ is Nash.
\item For all $v\in D$, $C_v\colon Aut_a(D)\to D\times Gl(n)$
is a Nash isomorphism onto its image.
\end{enumerate}
\end{Th}

	Theorem~\ref{cor} is now a corollary of Theorem~\ref{main}. This is
a consequence of the following result of K.~Diederich and
J.~E.~Forn\ae ss (\cite{DF}).

{\bf Theorem~1.4}
{\it Let $D\subset\subset\C^n$ be an algebraic domain. Then $Aut_a(D)=Aut(D)$.}

	{\bf Remark.} The proof of Theorem~1.4 makes use of the reflection
principle (see S.~Pin\v cuk, \cite{P}) and basic
methods of S.~Webster (\cite{W}) which are also fundamental in the
present paper. Due to the results of K.~Diederich and S.~Webster (\cite{DW})
and of S.~Webster (\cite{W}) on the continuation of automorphisms,
Theorem~\ref{main} can be applied to every situation where the
automorphisms can be $C^\infty$ extended to the boundary.

	K.~Diederich has informed the author that he and S.~Pin\v cuk recently
proved that, under natural non-degeneracy conditions on the boundary,
automorphisms of domains are always almost everywhere continuously
extendable. Applying this along with the reflection method, one would
expect $Aut(D)=Aut_a(D)$ for $D$ a semi-algebraic domain with $\partial D$
completely Levi-non-degenerate. Thus Theorem~\ref{main} applied to this
situation would show that $Aut(D)$ is an affine Nash  group acting
semi-algebraically on $D$.

	{\bf Acknowledgement.} On this occasion I would like to thank my teacher
A.~T.~Huckleberry for formulating the problem, for calling my attention
to the relevant literature and for numerous very useful discussions.

\section{Real semi-algebraic sets}\label{sa}

	Here we present some basic properties of semi-algebraic sets which
will be used in the proof of Theorem~\ref{main}.
For the proofs we refer to Benedetti-Risler \cite{BR}.

\begin{samepage}
\begin{Def}\label{def-s-a}
	A subset $V$ of \/ $\R^n$ is called {\bf semi-algebraic}
if it admits some representation of the form
$$V = \bigcup_{i=1}^s \bigcap_{j=1}^{r_i} V_{ij}$$
where, for each $i=1,\ldots,s$ and $j=1,\ldots,r_i$,
$V_{ij}$ is either $\{ x\in\R^n \mid P_{ij}(x)<0 \}$ or
$\{ x\in\R^n \mid P_{ij}(x)=0 \}$ for a real polynomial $P_{ij}$.
\end{Def}
\end{samepage}

	As a consequence of the definition it follows that finite unions and
intersections of semi-algebraic sets are always semi-algebraic. Moreover,
closures, boundaries, interiors (see Proposition~\ref{cl})
and connected components of semi-algebraic sets are semi-algebraic.
Further, the number of connected components is finite
(see Corollary~\ref{con}). Finally, any semi-algebraic set admits
a finite semi-algebraic stratification (see Definition~\ref{str} and
Proposition~\ref{strat}).

	The natural morphisms in the category of semi-algebraic set are {\bf
semi-algebraic maps}:

\begin{Def}\label{map}
	Let $X\subset\R^n$ and $Y\subset\R^n$ be semi-algebraic sets.
A map $f\colon X\to Y$ is called {\bf semi-algebraic} if the graph
of $f$ is a semi-algebraic set in $\R^{m+n}$.
\end{Def}

\begin{samepage}
\begin{Def}\label{str}
	A {\bf stratification} of a subset $E$ of $\R^n$ is a
partition $\{A_i\}_{i\in I}$ of $E$ such that
	\begin{enumerate}
	\item each $A_i$ (called a {\bf stratum}) is a real analytic
locally closed submanifold of $\R^n$;
	\item if $\overline{A_i} \cap A_j \ne \emptyset$,
then $\overline{A_i} \supset A_j$ and $\dim A_j < \dim A_i$ (frontier
condition).
	\end{enumerate}
	A stratification is said to be {\bf finite} if there is a finite
number of strata and to be {\bf semi-algebraic} if furthermore
each stratum is also a semi-algebraic set.
\end{Def}
\end{samepage}

\begin{Prop}\label{strat}
	Every semi-algebraic set $E\subset\R^n$ admits a semi-algebraic
stratification.
\end{Prop}

\begin{Cor}\label{con}
	Every semi-algebraic set has a finite number of connected components
and each such component is semi-algebraic.
\end{Cor}

\begin{Prop}\label{cl}
	Let $X$ be a semi-algebraic set in $\R^m$. Then the closure $\bar X$,
its interior $\mathop{X}\limits^0$, and its boundary $\partial X$
are semi-algebraic sets.
\end{Prop}

	Using Proposition~\ref{strat}, dimension of a semi-algebraic
set is defined to be the maximal dimension of its stratum.
This is independent of
the choice of a finite stratification.

\begin{Prop}\label{as}
	Let $Y\subset\R^m$ be a semi-algebraic set of $\dim Y\le k$.
Then it is contained in some real algebraic set $Z$ with ${\rm dim}Z\le k$.
\end{Prop}

	The following results on images and local triviality of semi-algebraic maps
will play an important role in the present paper.

\begin{Th}[Tarski-Seidenberg]\label{TS}
	Let $f\colon X \to Y$ be a semi-algebraic map. Then the image
$f(X)\subset Y$ is semi-algebraic set.
\end{Th}

	Further we need the Theorem on local triviality
(see Benedetti-Risler \cite{BR}, Theorem~2.7.1, p.~98).

\begin{samepage}
\begin{Th}\label{tr}
	Let $X$ and $Y$ be semi-algebraic sets and let $f\colon X\to Y$
be a continuous semi-algebraic map. Fix a finite semi-algebraic partition
of $X$, $\{X_1,\ldots,X_h\}$. Then there exists
	\begin{enumerate}
	\item a finite semi-algebraic stratification $\{Y_1,\ldots,Y_k\}$ of $Y$;
	\item a collection of semi-algebraic sets $\{F_1,\ldots,F_k\}$ and,
for every $i=1,\ldots,k$, a finite semi-algebraic partition
$\{ F_{i1},\ldots,F_{ir} \}$ of $F_i$ (typical fibres);
	\item a collection of semi-algebraic homeomorphisms
$$ g_i\colon f^{-1}(Y_i) \to Y_i\times F_i, \quad i=1,\ldots,k$$
such that
		\begin{enumerate}
		\item the diagram
$$\def\normalbaselines{\baselineskip20pt\lineskip3pt\lineskiplimit3pt}
\def\mapright#1{{\buildrel #1 \over \longrightarrow}}
\def\mapdown#1{\Big\downarrow \rlap{$\scriptstyle#1$}} \matrix{
f^{-1}(Y_i) & \mapright{g_i} & Y_i\times F_i \cr
\mapdown{f} &                &\mapdown{p}    \cr
Y_i         & \mapright{id}  & Y_i           \cr
}$$
is commutative ($p\colon Y_i\times F_i \to Y_i$ is the natural projection);
		\item for every $i=1,\ldots,k$ and every $h=1,\ldots,r$,
$$ g_i(f^{-1}(Y_i)\cap X_h) = Y_i \times F_{ih}. $$
		\end{enumerate}
	\end{enumerate}
\end{Th}
\end{samepage}

\begin{Rem}
	In the above setting one says that
$f$ is a trivial semi-algebraic map over $Y_i$ with
typical fibre $F_i$ and structure homeomorphisms $g_i$.
\end{Rem}

	Further, we shall use the following proposition from real algebraic
geometry (see \cite{BR}, Proposition~3.2.4).

\begin{Prop}\label{ri}
	Let $V\subset\R^n$ be a real algebraic variety. The set of singular
points is an algebraic set properly contained in $V$.
\end{Prop}

	The following Proposition provides a motivation for Definition~\ref{br}.

\begin{Prop}\label{bra}
	The Nash automorphisms $f\in Aut_a(D)$ are branches of algebraic maps.
\end{Prop}

	{\bf Proof.}
Let $f\in Aut_a(D)$ be a Nash automorphism. Then $f$ is holomorphic and
semi-algebraic. Every coordinate $f_j\colon D\to \C$ is also
holomorphic and, by Theorem~\ref{TS}, semi-algebraic.
By Proposition~\ref{as}, there
exist real algebraic sets $Z_j$ of (real) codimension $2$ with
$\Gamma_{f_j}\subset Z_j$. By Proposition~\ref{ri}, there exists a regular
point $(w_0,f(w_0))\in \Gamma_{f_j}$ where
$\Gamma_{f_j}$ is locally given
by two real polynomials $P_1(z,\bar z)$, $P_2(z,\bar z)$ or by a complex one
$P(z,\bar z):=P_1(z,\bar z)+iP_2(z,\bar z)$, such that $dP\ne 0$.
The latter property implies that
either $\partial P\ne 0$ or $\bar\partial P\ne 0$.
We can assume that $\partial P\ne 0$, otherwise $P$ can be replaced with $\bar
P$.

	Let $P(z,\bar z)=P'(z)+P''(z,\bar z)$ be a decomposition of $P$ in
the holomorphic part $P'$ and the remainder $P''$ which
consists only of terms with non-trivial
powers of $\bar z$. The part $P'$ is not zero because
$dP'=\partial P\ne 0$.
We wish to prove that $\Gamma_{f_j}$ is locally defined by
the holomorphic polynomial $P'(z)$. The identity
$$P(w,f_j(w),\bar w, \overline{f_j(w)}) \equiv 0$$
for all $w$ near $w_0$ implies in particular the vanishing of
the Taylor coefficients of $w^k$ for all multi-indices $k$.
But these coefficients are just
$${1\over |k|!}{\partial^{|k|} P'(w,f_j(w)) \over \partial w^k}$$
and their vanishing means the vanishing of $P'(z,f_j(z))$.

	Thus the graph $\Gamma_{f_j}$ is locally defined by
a holomorphic polynomial $P_j$
and the polynomials $P_1,\ldots,P_n$ define the $n$-dimensional
algebraic variety required in the definition of a branch of an algebraic map.
\nopagebreak\par\hfill {\bf Q. E. D.}

	We shall also make frequent use of Chevalley's theorem on
constructible sets (see Mumford \cite{M}, p.~72).

\begin{Def}\label{cons}
	A subset $A\in\C^n$ is called {\bf constructible} if it is a finite
union of locally closed complex algebraic subvarieties.
\end{Def}

\begin{Rem}
	Constructible sets are semi-algebraic.
\end{Rem}

\begin{Th}[Chevalley]\label{Chev}
	Let $X$ and $Y$ be affine varieties and $f\colon X\to Y$ any morphism.
Then $f$ maps constructible sets in $X$ to constructible sets in $Y$.
\end{Th}

\section{A scheme of the proof}

	The proof of Theorem~\ref{main} can be divided in two steps.
The essential ingredient for the first step is the method of S.~Webster
(see \cite{W}) based on the reflection principle
(see S.~Pin\v cuk, \cite{P}). We use it to construct
an appropriate {\it family of graphs} of automorphisms from $Aut_a(D)$.
In fact we construct a constructible family $F$ which fibres contain
automorphisms from $Aut_a(D)$.
This is carried out in sections~\ref{wm} and \ref{we}.

  In the second step we show using the family constructed in
the first step that the set $C_v(Aut_a(D))$ is Nash.
Taking a neighborhood where $C_v(Aut_a(D))$ is not empty and closed and
taking its pullback in $Aut(D)$ we obtain a neighborhood where $Aut_a(D)$
is closed, which implies statement~1. in Theorem~\ref{main}.
In fact we prove that the {\it exact}
family of graphs
$$\Gamma:=\{(C_v(f),w,f(w) \mid f\in Aut_a(D), \, w\in\C^n \}$$
is Nash. If we identify $Aut_a(D)$ with $C_v(Aut_a(D))$,
$\Gamma$ is the graph of the action $Aut_a(D)\times D\to D$. This
proves statement~3. To obtain statement~2., we observe that group operation
can be defined {\em semi-algebraically} in terms of $\Gamma$.
Here we use the theory of semi-algebraic sets and their morphisms.

	\section{Reflection principle}\label{wm}

	Let $D\subset\subset\C^n$ be a semi-algebraic Levi-non-degenerate domain.
By Proposition~\ref{cl}, the boundary $\partial D$ is
semi-algebraic. By Proposition~\ref{as}, there exists a real algebraic
set $H$ of dimension $2n-1$ which contains $\partial D$. Let
$H_i$ be irreducible components of $H$ of dimension $2n-1$.
By Definition~\ref{L}, the Levi form of some component of $H$, let say
of $H_1$, is not everywhere degenerate.

	To every irreducible hypersurface $H_i$ we associate a real Zariski open
set $U_i\subset\C^n$ and a real polynomial $r_i(z,\bar z)$ with
$H_i\cap U_i= \{r_i=0\} \cap U_i$ and $dr\ne 0$ on $U_i$.
By Proposition~\ref{ri}, such $U_i$'s and $r_i$'s exist.

	Let $f\in Aut_a(D)$ be any fixed map which is, by Proposition~\ref{bra}
a branch of an
algebraic map and $V\subset\C^n\times\C^n$ be the corresponding
$n$-dimensional algebraic subvariety which contains the graph of $f$.
We wish to extend $f$ in a neighborhood of a boundary point $x\in H_1$.
Outside a proper complex algebraic subvariety $V'\subset V$ the variety
$V$ defines $f$ and $f^{-1}$ as possibly multiple-valued algebraic maps.
Since $\dim_{\C} V=n$, $\dim_{\C} V'\le n-1$. Since $\dim_{\R} H_1=2n-1$,
there exists a point $x\in H_1\cap U_1$ and
a neighborhood $U\subset U_1$ of $x$ such
that $V$ is a trivial covering over $U$
and all sections of this covering define biholomorphic maps onto their images.
One of these maps coincides with $f$ over $D\cap U$. This map yields the
desired extension of $f$. Since $f$ is an automorphism of $D$,
it maps $H_1\cap U$ into $\partial D$.

{\bf Notation.} Let $i=i(f)$ be such that $f(H_1\cap U)\subset H_i$.
Let $w_0\in H_1\cap U$ be
an arbitrary point such that for $w'_0=f(w_0)$ one has $dr_i(w'_0)\ne 0$.

We use the notation
$$ r:=r_1, r':=r_i, H:=H_1 \hbox{, and } H':=H_i.$$
Then $f(H\cap U)\subset H'$ and we have a relation
\begin{equation}\label{pz}
	r'(f(z),\bar f(\bar z)) = g(z,\bar z) r(z,\bar z),
\end{equation}
where $g(z,\bar z)$ is real analytic.

	Let $z=x+iy$, where $x$ and $y$ are real coordinate vectors. Since
the functions in (\ref{pz}) are given by power series in $(x,y)$,
they are still defined for complex vectors $x$ and $y$ near $w_0$.
This is equivalent to varying $z$ and $\bar z$
independently. The relation (\ref{pz}) persists:
\begin{equation}\label{pzw}
	r'(f(z),\bar f(\bar w)) = g(z,\bar w) r(z,\bar w).
\end{equation}

	Now we consider the spaces $Z:=\C^n$, $Z':=\C^n$, $W:=\C^n$
and $W':=\C^n$ and define
the complexifications ${\cal H}\subset Z\times\bar W$ and
${\cal H'}\subset Z'\times\bar W'$ by
\begin{equation}\label{wz}
	r(z, \bar w) = 0, \quad r'(z', \bar w') = 0.
\end{equation}

	The so-called {\em Segre complex varieties} associated to
the points $w\in W$ are defined by
$$ Q_w = \{ z\in Z \mid r(z,\bar w)=0 \}. $$
Since $g(z,\bar w)$ is holomorphic for $z$ and $w$ near $w_0$,
we see from (\ref{pzw}) that the map $f^{\C}:=f\times\bar f$ takes
$\{z\} \times Q_z$ into $\{z'\} \times Q'_{z'}$,
where $z'=f(z)$. Hence, the family of complex hypersurfaces
$\{z\} \times Q_z$ is invariantly related to $H$.

	Since $r(z,\bar z)$ is real, we have
$$ r(z,\bar w) = \bar r(\bar w, z) = \overline{r(w,\bar z)}, $$
so that $z\in Q_w \iff w\in Q_z$. Also $z\in Q_z \iff z\in H$.
Since $r$ is real and

\newpage
\noindent $dr=\partial r + \bar\partial r$
does not vanish at $w_0$ (resp. $r'$ is real and
$dr'=\partial r' + \bar\partial r'$ does not vanish at $w'_0$), we have
\begin{equation}\label{dr}
	\partial r(z,\bar w)\ne 0, \quad \partial r'(z',\bar w')\ne 0
\end{equation}
for $z$ and $w$ near $w_0$ and for $z'$ and $w'$ near $w'_0$.

	The relation (\ref{dr}) implies that $Q_w$ is non-singular
in $z$ if $(z,\bar w)$ is near $(w_0,\bar w_0)$.
Let $\pi_z(\bar w)$ denote the complex tangent
space $T_z Q_w$ as an element in the grassmanian $G_{n,n-1}$.
It follows that $\pi_z$ is an
antiholomorphic map from $\{z\}\times Q_z$ to $G_{n,n-1}$.
S.~Webster proves the following fact (\cite{W}, p.~55, Lemma~1.1):

\begin{Lemma}\label{inv}
	The antiholomorphic map $\pi_z(\bar w)$ is locally
invertible near the points of $H$ where the Levi form is non-degenerate.
\end{Lemma}

	Since the Levi form of $U\cap H$ is non-degenerate, and
a biholomorphic map preserves this property, the Levi form of $H'$ is
non-degenerate around $w'_0=f(w_0)\in H'$.
By Lemma~\ref{inv}, $\pi_z$ and $\pi'_{z'}$ are
locally invertible for $(z,\bar w)$ near $(w_0,\bar w_0)$ and
$(z',\bar w')$ near $(w'_0,\bar w'_0)$.

	The first step of Webster's method is to describe
the map $f$ between $Q_z$ and $Q'_{z'}$. We have seen
that $f$ takes $Q_w$ into $Q'_{w'}$, $w'=f(w)$. We take $w\in Q_z$.
Then all $Q_w$'s pass through the point $z$ and all $Q'_{w'}$'s through $z'$.
Therefore the differential $f_{*z}\in Gl(n)$ takes
$T_zQ_w$ into $T_{z'}Q'_{w'}$,
i.e. $\pi_z(\bar w)$ into $\pi_{z'}(\bar w')$. These considerations mean
that the restriction $f|_{Q_z}$ can be decomposed as follows (\cite{W}, p.~56):
\begin{equation}\label{dec}
	w \stackrel{\vphantom{\pi_{z'}^{-1}}\pi_z}{\longmapsto} T_zQ_w
	\stackrel{\vphantom{\pi_{z'}^{-1}}J_{n-1}(f_{*z})}{\longmapsto}
	T_{z'}Q'_{w'} \stackrel{\pi_{z'}^{-1}}{\longmapsto} w',
\end{equation}
where the second map $J_{n-1}(f_{*z})\colon G_{n,n-1}\to G'_{n,n-1}$
is the natural map between grassmanians which is induced by the differential
$f_{*z}\colon T_zZ \to T_{z'}Z'$.

	The decomposition (\ref{dec}) implies that for $w\in Q_z$ and $w'=f(w)$,
$z'=f(z)$ and $q=f_{*z}\in Gl(n)$ the following relation is satisfied:
\begin{equation}\label{wqw}
	\pi'_{z'}(\bar w') = J_{n-1}(q)(\pi_z(\bar w)).
\end{equation}

	This relation expresses the restriction $\phi:=f|_{Q_z}$
in terms of parameters $(z,z',q)\in Z\times Z'\times Gl(n)$.
This is of great importance for our parametrization of the graphs of
elements of $Aut_a(D)$. Hence we underline this fact by introducing
the notation
\begin{equation}\label{phi}
	\phi(\bar z,\bar z',\bar q)\colon Q_z \to Q'_{z'}
\end{equation}
for $\phi:=f|_{Q_z}$. We write the conjugate variables for arguments
of $\phi$ in order to emphasize that $\phi$ depends holomorphically on them.

	The idea of the second step is to express the map $f$ in terms of
restrictions $f|_{Q_z}$. This is done separately for $n=2$ and
$n\ge 3$.

	\subsection{The case $n\ge 3$}

	Let $z_0\in U\cap Q_{w_0}$ be any point.
Lemma~\ref{inv} gives points $v_1,\ldots,v_n\in Q_{z_0}$ near
$w_0$ such that all $Q_{v_j}$ are non-singular and transverse in $z_0$ and
the algebraic curve $\gamma$, defined by
\begin{equation}\label{zv0}
	r(z,\bar v_1) = \cdots = r(z,\bar v_{n-1}) = 0,
\end{equation}
is transverse to $Q_{w_0}$.
Every $w$ near $w_0$ lies in some $Q_z,z\in\gamma$, namely for
$z\in \gamma\cap Q_w$. Therefore, to describe $f(w)$ we need only to
consider restrictions $f|_{Q_z},z\in \gamma$.

	Given the values $z'=f(z)$ and differentials $q=f_{*z}\in Gl(n)$,
the map $f|_{Q_z}$ is determined by (\ref{wqw}).
Since $\gamma\subset Q_{v_1}$, the values $z'=f(z)$ along
$\gamma$ are determined, in turn, due to (\ref{wqw}) by parameters
$(v_1,v'_1,l_1)=(v_1,f(v_1),f_{*v_1})$.
Namely, we use the map $\phi$ in (\ref{phi}) and set
\begin{equation}\label{zz'}
	z'=\phi(\bar v_1,\bar v'_1,\bar l_1)(z).
\end{equation}

	Further, the differentials $f_{*z}$ along $\gamma$ can be expressed
in terms of parameters $(v_1,v'_1,l_1)=(v_1,f(v_1),f_{*v_1})$.
Consider the differential 1-forms
$$\theta_\alpha=\partial r(z,\bar v_\alpha).$$
They define a frame in the cotangent spaces.
Let
$$\{Y_j=Y_j(z,\bar v_1,\ldots,\bar v_{n-1},\bar w),j=1,\ldots,n\}$$
be the dual vector field frame.
This frame has rational coefficients in the
variables $(z,\bar v_1,\ldots,\bar v_{n-1},\bar w)$ and satisfies
the conditions
\begin{equation}\label{Y}
\left.\matrix{
	\hbox{$Y_1$ is transverse to $Q_{v_1}$ and tangent to $Q_{v_2}$,} \cr
	\hbox{$Y_2$ is transverse to $Q_{v_2}$ and tangent to $Q_{v_1}$,} \cr
	\hbox{$Y_3,\ldots,Y_n$ are tangent to $Q_{v_1} \cap Q_{v_2}$.}\cr
}\right\rbrace
\end{equation}

	Similar differential $1$-forms $\theta'_\alpha$
and frame vector fields $Y'_j$ are constructed for $H'$.
Relative to these two frame fields
$$ f_{*z}Y_l = \sum q_{lj}Y'_j, $$
where
\begin{equation}\label{mat}
   [q_{lj}]= \left[ \matrix{ q_{11}     & 0          & 0               \cr
			     0          & q_{22}     & 0               \cr
			     q_{1\beta} & q_{2\beta} & q_{\alpha\beta} \cr } \right],
\end{equation}
$\alpha,\beta=3,\ldots,n$.

	The functions $q_{11}$, $q_{1\beta}$, $q_{\alpha\beta}$ are determined by
values of $f$ along $Q_{v_1}$, i.e. by $\phi(v_1,v'_1,l_1)$
(where $v'_1=f(v_1)$ and $l_1=f_{*v_1}$).
Similarly, $q_{22}$ and $q_{2\beta}$ are determined by
$\phi(v_2,v'_2,l_2)$. These dependencies can be expressed by
relations
\begin{equation}\label{q}
\left.\matrix{
	q_{11}         & = &\theta'_1(\phi_{*z}(v_1,v'_1,l_1)Y_1),             \cr
	q_{1\beta}     & = &\theta'_1(\phi_{*z}(v_1,v'_1,l_1)Y_\beta),         \cr
	q_{\alpha\beta}& = &\theta'_\alpha(\phi_{*z}(v_1,v'_1,l_1)Y_\beta),    \cr
	q_{22}         & = &\theta'_2(\phi_{*z}(v_2,v'_2,l_2)Y_2),             \cr
	q_{2\beta}     & = &\theta'_2(\phi_{*z}(v_2,v'_2,l_2)Y_\beta),         \cr
	q_{12}         & = & q_{21} = q_{\alpha 1} = q_{\alpha 2}   = 0.       \cr
}\right\rbrace.
\end{equation}
	Thus, the map $f$ is completely determined by parameters
$v_j\in V_j:=\C^n$, $v'_j\in V'_j:=\C^n$ and $l_i\in L_j:=Gl(n)$,
$j=1,\ldots,n$.

	\subsection{The case $n=2$}

	In case $n=2$ there are no frames with
properties~(\ref{Y}) and another construction (\cite{W},~p.~58) is needed.
By Lemma~\ref{inv}, two points $\zeta_1,\zeta_2\in Q_{w_0}$ can be chosen
such that $Q_{\zeta_1}$ and $Q_{\zeta_2}$ are non-singular and transverse in
$w_0$. Then choose $v_1\in Q_{\zeta_1}$ and $v_2\in Q_{\zeta_2}$ such that
each $Q_{v_j}$ is non-singular in $\zeta_j$ and transverse to
$Q_{w_0}$ there. Now fix $v_1$ and $v_2$ and let $z_1$ and $z_2$ move
along $Q_{v_1}$ and $Q_{v_2}$ respectively. For $z_1$ and $z_2$ near
$\zeta_1$ and $\zeta_2$, it follows that $Q_{z_1}$ and $Q_{z_2}$ are still
transverse near $w_0$ and intersect each other in a single point $w$
there. Conversely, for given $w$ near $w_0$, $Q_w$ intersects
each $Q_{v_j}$ transversely in a point $z_j$ near $\zeta_j$.
In this way a local biholomorphic correspondence between $w\in W(:=\C^n)$ and
$(z_1,z_2)\in Q_{v_1}\times Q_{v_2}$ is obtained. It is defined by relations
\begin{equation}\label{etaw}
	r(w,\bar{z}_j) = 0 \quad (\iff z_j\in Q_w), \quad j=1,2;
\end{equation}
\begin{equation}\label{etav}
	r(v_j,\bar{z}_j) = 0 \quad (\iff z_j\in Q_{v_j}).
\end{equation}

	Further, set $w'_0:=f(w_0)$, $\zeta'_j:=f(\zeta_j)$, $v'_j:=f(v_j)$.
All transverse properties are preserved by the biholomorphic map $f$. Again,
one obtains a local biholomorphic correspondence between $w'\in W'$ and
$(z'_1,z'_2)\in Q_{v'_1}\times Q_{v'_2}$, which is defined by
\begin{equation}\label{etaw'}
	r'(w',\bar{z}'_j) = 0 \quad (\iff w'\in Q'_{z'_j});
\end{equation}
\begin{equation}\label{etav'}
	r'(v'_j,\bar{z}'_j) = 0 \quad (\iff z'_j\in Q'_{v'_j}).
\end{equation}

	Since the Segre varieties $Q_z$ are invariant with respect to $f$,
if $z'_j:=f(z_j)$, one obtains the corresponding point $w'=f(w)$.
Thus, $f$ can be decomposed in the following way:

\begin{equation}\label{dec1}
	W \longrightarrow Q_{v_1}\times Q_{v_2} \stackrel{f\times f}
	{\longrightarrow} Q'_{v'_1}\times Q'_{v'_2} \longrightarrow W'.
\end{equation}

	The middle map here is in fact $f|_{Q_{v_1}} \times f|_{Q_{v_2}}$ which
is equal to
$\phi(\bar v_1,\bar v'_1,\bar l_1) \times \phi(\bar v_2,\bar v'_2,\bar l_2)$,
where $\phi$ is the map (\ref{phi}) and $l_j:=f_{*v_j}$. In other words
we have a relation between $z_j$ and $z'_j$:

\begin{equation}\label{zz'1}
	z'_j = \phi(\bar v_j,\bar v'_j,\bar l_j)(z_j)
\end{equation}

Thus, $f$ is completely determined by parameters
$v_j\in V_j:=\C^n$, $v'_j\in V'_j:=\C^n$ and $l_i\in L_j:=Gl(n)$, $j=1,2$.

\subsection{Reflection principle with parameters}\label{we}

	The local construction recalled in previous paragraph
is in fact global because of its algebraic nature. The map $f$
was locally expressed in terms of parameters
$(v,v',l)=(v,f(v),f_{*v})\in P$, where $v:=(v_1,\ldots,v_n)$ and
$$P:=\prod_{j=1}^n (V_j\times V'_j\times L_j).$$
Using the same algebraic relations
globally, we shall obtain a constructible family
$F\subset P\times W\times W'$
such that the graph $\Gamma_f$ is an open subset of the closure of the fibre
$\overline{F_p}$ for generic $v$ and $p=(v,f(v),f_{*v})$
(this will be made precise below).

	We start with construction of a constructible family for the map
$$\phi(z,z',q)\colon Q_z \to Q'_{z'}$$
in (\ref{phi}). For this we consider the
constructible subset
$$\Phi \subset \overline{Gl(n)} \times \bar Z \times W \times \bar Z' \times
W'$$
defined by relations (\ref{wz}), (\ref{dr}), and (\ref{wqw}).
The relations (\ref{dr}) provide the existence of
$\pi_z(\bar w)$ and $\pi'_{z'}(\bar w')$ respectively. Furthermore, we
have seen that
$$(\overline{f_{*z}},\bar z,w,\overline{f(z)},f(w))\in \Phi$$ for
$(z,\bar w)\in \cal H$ near $(w_0,\bar w_0)$, and $\pi_z$ and $\pi_{z'}$ are
local invertible there (Lemma~\ref{inv}). To provide this local
invertibility ``globally'',
we assume, changing if necessary to a smaller constructible subset of $\Phi$,
that
\begin{equation}\label{dpi}
	\det {\partial \pi_z(\bar w) \over \partial \bar w} \ne 0,
\quad   \det {\partial \pi'_{z'}(\bar w') \over \partial \bar w'} \ne 0.
\end{equation}

	The set $\Phi$ defines now a family of possibly multiple-valued maps
$$\phi(\bar z,\bar z',\bar q)\colon Q_z \to Q'_{z'}.$$
Consider the complexifications ${\cal H}\subset Z\times \bar W$ and
${\cal H'}\subset Z'\times \bar W'$. Then $\Phi$ is a subset in
$\overline{Gl(n)} \times \bar{\cal H} \times  \bar{\cal H}'$.

\begin{Lemma}\label{f1}
	The projection
$\delta\colon \Phi \to \overline{Gl(n)} \times \bar Z' \times \bar {\cal H}$
has finite fibres and is locally biholomorphic.
\end{Lemma}

	{\bf Proof.} We need to prove that
$w'\in\delta^{-1}(\bar q,\bar z',\bar z,w)$
depends locally holomorphically on $(\bar q,\bar z',\bar z,w)$.
For this it is enough to observe, that, by (\ref{dpi}),
$\pi_{z'}(\bar w)$ is locally invertible and, by (\ref{wqw}),
$\bar w'=(\pi'_{z'})^{-1}(J_{n-1}(q)(\pi_z(\bar w)))$. Since the above
fibres are constructible, they are finite.
\par\hfill {\bf Q. E. D.}

	In the following let $\phi$ denote the multiple-valued map defined
by $\Phi$. Since every value of $\phi$ is, by Lemma~\ref{f1},
locally holomorphic in $w$,
we can discuss its differential $\phi_*=\phi(\bar q,\bar z',\bar z,w)_*$
which is also possibly multiple-valued.

	For the construction of the required family we need to consider
auxiliary parameter spaces $A:=Z\times Z'\times Gl(n)$
for $n>2$ and $A:=(Z\times Z')^2$ for $n=2$.
Let $F\subset \bar A \times P \times W\times W'$ be
the constructible subset defined by relations (\ref{wz}),
(\ref{wqw}), (\ref{zv0}), (\ref{zz'}) and (\ref{q})
in case $n>2$ and by (\ref{etaw}), (\ref{etav}), (\ref{etaw'}), (\ref{etav'})
and (\ref{zz'1}) in case $n=2$.

	Passing if necessary to a constructible subset, we can require that
in case $n>3$ all $Q_{v_j}$'s and $Q_w$ are transverse in $z$ and
all $Q'_{v'_j}$'s and $Q'_{w'}$ are transverse in $z'$.
In case $n=2$ we require that each $Q_{v_j}$ is transverse to $Q_w$ in $z_j$
and $Q_{z_j}$'s are transverse in $w$ and, similarly,
each $Q'_{v'_j}$ is transverse to $Q'_{w'}$ in $z'_j$
and $Q'_{z'_j}$'s are transverse in $w'$.

	Further, by Theorem~\ref{Chev} of Chevalley,
the projection $\pi(F)$ of $F$ on $P\times W\times W'$ is also constructible.
We don't have in general a local biholomorphic
property as in Lemma~\ref{f1} for $\pi(F)$,
but we still can prove the finiteness:

\begin{Lemma}\label{f2}
	The projection $\sigma\colon \pi(F) \to P\times W$ has finite fibres.
\end{Lemma}

	{\bf Proof.} Let fix $(p,w)\in P\times W$.
Let $(p,w,w')\in F$ be any point. By the construction of $\pi(F)$,
there exist points $a\in A$ such that $(\bar a,p,w,w')\in F$.

	{\em Case $n>2$} Let $a=(z,z',q)$. We constructed $F$ such that
$Q_{v_1},\ldots,Q_{v_{n-1}}$ and $Q_w$ are transversal in $z$. Then, by
(\ref{wz}) and (\ref{zv0}), the set of possible $z\in Z$ is discrete and
therefore finite. Further, by (\ref{zz'}) and (\ref{q}), only finitely
many $z'$'s and $q$'s are possible. Here we use Lemma~\ref{f1}.
Now $w'\in W'$ is determined by (\ref{wqw}), which implies finiteness of
the set of $w'$'s.

	{\em Case $n=2$} Let $a=(z_1,z'_1,z_2,z'_2)$. By definition of $F_2$,
$Q_w$ and $Q_{v_j}$ are transverse in $z_j$. Then there are only finitely
many possible intersections $z_j$. By (\ref{zz'1}),
the number of possible $z'_j$'s is also finite. Finally, since $Q'_{z'_j}$ are
transverse in $w'$, the number of possible $w'$'s is also finite.
\par\hfill {\bf Q. E. D.}

	Now let $f\in Aut_a(D)$ be fixed and $H=H_1$, $H'=H_i$ for $i=i(f)$
(we defined $i(f)$ by the condition $f(H_1\cap U)\subset H_i$ for some open
$U\subset\C^n$ with $H_1\cap U\ne\emptyset$).
By our construction,
$(v,f(v),f_{*v},w,f(w))\in \pi(F)$ for all $(v,w)$ in some open subset
$U_1\subset V\times W$. Here we wish to point out that the family
$\pi(F)$ depends on the
index $i=i(f)$. To include all automorphisms $f\in Aut_a(D)$, we just
consider the finite union of $\pi(F)_i\subset P\times W\times W'$ for
all possible $i=i(f)$, $f\in Aut_a(D)$ and denote it again by $F$.

	The set $F$ is constructible, i.e. a finite union
of locally closed algebraic subvarieties. It follows that the set
$$E(f):= \{(v,w)\in D^{n+1} \mid (v,f(v),f_{*v},w,f(w))\notin F \} $$
is analytically constructible, i.e. a finite union of locally closed
analytic subvarieties. If $(v,w)$ is outside the closure of $E(f)$, we have
$(p,w,f(w))\in F$ for $p=(v,f(v),f_{*v})\in P$. This means that
the graph $\Gamma_f:=\{(w,f(w)) \mid w\in D \}$ lies in the closure of
the fibre $F_p$.

	Now we wish to prove main result of this section.

\begin{samepage}
\begin{Prop}\label{par}
	Let $P$, $W$ and $W'$ be as above. There exist constructible subsets
$F\subset P\times W\times W'$ and $E\subset P\times V\times W$ such that
	\begin{enumerate}
	\item The projection $\sigma\colon F\to P\times W$ has finite fibres;
	\item For every fixed $f\in Aut_a(D)$ there exists a proper subset
$E(f)\subset D^{n+1}$, such that for all
$(v,w)\in (D^{n+1})\backslash E(f)$ and $p=(v,f(v),f_{*v})\in P$
one has $(p,w,f(w))\in F$, the graph $\Gamma_f$ is a subset
of the closure of the fibre $F_p$ and $E(f)\subset E_p$;
	\item $E_p$ is of complex codimension at least $1$ in $V\times W$.
	\end{enumerate}
\end{Prop}
\end{samepage}

	We need the following lemma.

\begin{Lemma}\label{clos}
	Let $A$, $B$ and $C\subset A\times B$ be constructible subsets of
arbitrary algebraic varieties. Then the
{\bf fibrewise} closure of $C$ in $A\times B$,
i.e. the union of closures of the fibres
$C_a := (\{a\} \times B) \cap C$, $a\in A$ is constructible.
\end{Lemma}

	The proof is based on the following fact
(see Mumford, \cite{M}, Corollary~1, p.71). Recall that a morphism is
{\it dominating} if its image is dense.

\begin{Prop}\label{mor}
	Let $X$ and $Y$ be two complex algebraic varieties,
$f\colon X\to Y$ be a dominating morphism and $r=\dim X-\dim Y$.
Then there is a nonempty open set $U\subset Y$ such that, for all $y\in U$,
$f^{-1}(y)$ is a nonempty ``pure'' $r$-dimensional set, i.e. all its
components have dimension $r$.
\end{Prop}

	{\bf Proof of Lemma~\ref{clos}.}
	We first observe that given two constructible subsets
$C_1,C_2\subset A\times B$ which have constructible fibrewise closures,
the union $C=C_1\cup C_2$ has also this property. Changing to locally
closed irreducible components, we can assume that $A$, $B$ and $C$ are
irreducible algebraic varieties.

	Now we prove the statement by induction on dimension of $A$.
In case $\dim A=0$ the fibrewise closure of $C$ is just the closure of $C$
which is constructible.

	Let $\pi\colon \bar C\to A$ denote the projection of the closure
$\bar C$ on $A$. We can assume $\pi$
to be dominant, otherwise $A$ is replaced by the closure of $\pi(\bar C)$
which has a smaller dimension. Now we apply Proposition~\ref{mor} to
the projection $\pi$ and obtain an open subset $U\subset A$, such that
the fibre's over $U$ have pure dimension $\dim C-\dim A$. We have a
partition $A=U\cup (A\backslash U)$ of $A$ and the corresponding
partition $C=(C_1\cup C_2)$
($C_1:=C\cap (U\times B)$, $C_2:=C\cap ((A\backslash U)\times B)$).
By the above observation, it is enough to
prove the statement for $C_1$ and $C_2$ separately.
The statement for $C_2$ follows by induction, because
$\dim (A\backslash U)< \dim A$. Therefore we can assume $A=U$.

	Now we consider the irreducible components $C_i$ of
$\bar C \backslash C$, $\dim C_i<\dim C$.
If $S_i:=\overline{\pi(C_i)}\ne A$ for some $i$, then we replace $A$ by
$A\backslash S_i$ and correspondingly $C$ by $C\cap\pi^{-1}(A\backslash S_i)$.
Thus we may assume
that $\pi\colon C_i\to A$ is dominaiting for all $i$.
Then we can apply Proposition~\ref{mor} to
every $C_i$ and obtain a number of open sets $U_i\subset A$.
Let $U$ be the intersection of all $U_i$'s. Since $A$ is irreducible,
$U$ is not empty. Again, proceeding by induction, we can reduce the
statement to the case $A=U$. But in this case the fibres of $\bar C$ are
of pure dimension $\dim C-\dim A$ and the fibres of
$\bar C \backslash C$ have smaller dimension. This implies that
the fibrewise closure of $C$ coincides with the usual closure $\bar C$
which is constructible.
\par\hfill {\bf Q. E. D.}

	{\bf Proof of Proposition~\ref{par}.} Statement~1 follows from Lemma~\ref{f2}.
It follows from the local Webster's construction (section~\ref{wm}) that
$(p,w,f(w))\in F$ ($p=(v,f(v),f_{*v})$) for all $(v,w)$ in an open subset
$U\subset D^{n+1}$. This means that $U$ lies in the complement of the
`` exceptional set'' $E(f)$. Let
$\Omega(f):=D^{n+1}\backslash E(f)$, i.e.
$$\Omega(f)= \{(v,w)\in D^{n+1} \mid (v,f(v),f_{*v},w,f(w))\in F \}. $$

	The set $E$ must be globally defined independently of any automorphism
$f\in Aut_a(D)$. For this it is necessary to define $\Omega(f)$ in another
way. Changing if necessary to a constructible subset of $F$,
we may assume that the projection
$\sigma_p \colon F_p\to W$ is locally biholomorphic and $\Omega(f)$ still
contains an open subset $U\subset D^{n+1}$. Then the
differentials $\partial w' \over \partial w$ are certainly defined.
We now define the family $F'\subset F\times Gl(n)$
of differentials by adding values of
$\partial w' \over \partial w$:
$$F':= \{ (p,w,w',q)\in F\times Gl(n) \mid q={
\partial w' \over \partial w} \}$$
This is a constructible set and we have
$(v,f(v),f_{*v},w,f(w),f_{*w})\in F'$ for $(v,w)$
in some open set $U\subset D^{n+1}$.
Now we write the definition of $\Omega(f)$ in the form:
$$\Omega(f):= \{ (v,w) \in V\times W
\mid \hbox{ for } v'=f(v),l=f_{*v},w'=f(w): $$
$$ (v,v',l,w,w')\in F \}. $$
Now we define a set $\Omega$ which contains $\Omega(f)$ for all $f\in
Aut_a(D)$:
\begin{samepage}
$$\Omega:= \{ (p,v,w,v',w',l) \in P\times V\times W\times V'\times W'\times
Gl(n)
\mid $$
$$\forall j: (p,v_j,v'_j,l_j)\in F' \land (p,w,w')\in F \land
(v,v',l,w,w')\in F \}. $$
\end{samepage}
For $f\in Aut_a(D)$, $(v,w)\in U$, and $p=(v,f(v),f_*v)$,
we have
$$(v,w,f(v),f(w),f_{*v})\in \Omega_p.$$
Let $\Omega'$ be the fibrewise closure of $\Omega$, i.e. the union of all
closures
of $\Omega_p$, $p\in P$. By Lemma~\ref{clos}, $\Omega'$ is constructible.
Finally, we define
$E\subset P\times V\times W$ to be the projection of $\Omega'\backslash \Omega$
on $P\times V\times W$. By Theorem~\ref{Chev} of Chevalley, $E$ is
constructible.

	We now wish to prove that every fibre $E_p$ is of codimension at
least $1$. For this we note that for $p\in P$ fixed the projection of
$\Omega_p$ on $V\times W$ has finite fibres,
i.e. $\dim \Omega_p\le \dim (V\times W)$.
This implies $\dim (\Omega'\backslash \Omega)_p < \dim (V\times W)$
and $E_p$ is of codimension at least $1$ as required in statement~3.

	We take now any $(v,w)\in D^{n+1}$ outside $E(f)$ and set $p=(v,f(v),f_{*v})$.
For the proof of statement~2, consider $f\in Aut_a(D)$, take
$(v,w)\in D^{n+1}\backslash E(f)$, and set $p=(v,f(v),f_*v)$.
We shall prove that $E(f)\subset E_p$.
Let $(v_0,w_0)\in E(f)$ be any point. If $(v,w)\in U$, we have
$$(v,f(v),f_{*v},w,f(w),f_{*w})\in F',$$
which implies
$$(p, v, w, f(v), f(w), f_{*v}) \in \Omega.$$
Here $(v,w)\in U$ is arbitrary. Since $f$ is holomorphic,
we have this property globally for all $(v,w)\in D^{n+1}$ if we replace
$\Omega$ with its fibrewise closure $\Omega'$. In particular, we have
$$(v_0,w_0,f(v_0),f(w_0),f_{*v_0}) \in \Omega'_p.$$
Since $(v_0,w_0)\in E(f)$, the point
$$(v_0,f(v_0),f_{*v_0},w_0,f(w_0))$$
does not lie in $F$. This implies that
$$(v_0,w_0,f(v_0),f(w_0),f_{*v_0})$$
does not lie in $\Omega_p$ and then
it is in $\overline{\Omega_p}\backslash \Omega_p$. This means $(v_0,w_0)\in
E_p$,
which is required. The proof of Proposition~\ref{par} is finished.
\par\hfill {\bf Q. E. D.}

	\section{The choice of parameters}

	In the previous section we proved the existence of a constructible
algebraic family $F\subset P\times W\times W'$ with the
property that
for all $f\in Aut_a(D)$ there exists a point $p\in P$ such that
\begin{equation}\label{subs}
	\Gamma_f\subset \overline{F_p}.
\end{equation}
The goal of this section is to choose for every $f$ appropriate parameter $p$
with this property and obtain a map $\imath$ from $Aut_a(D)$ in the
corresponding
parameter space $P$. The first idea is to take some generic
$v\in V(:=V_1\times\cdots\times V_n=\C^{n^2})$ and to define
$p=(v,f(v),f_{*v})$. If $(v,w)\notin E(f)$, Proposition~\ref{par} yields the
required property (\ref{subs}). However, if we wish to define a global map
$Aut_a(D)\to P$, the condition $(v,w)\notin E(f)$ must be satisfied for
all $f\in Aut_a(D)$. Unfortunately, this is not true in general. It is
therefore necessary to take sufficiently many points
$(v_\mu,w_\mu) \in V\times W$ instead of
one $(v,w)$, such that $(v_\mu,w_\mu) \notin E(f)$ is always true at least
for one $\mu$. In fact, we prove the following Proposition.

\begin{samepage}
\begin{Prop}\label{par1}
	There exists a natural number $N$, a constructible subset
$F\subset P^N\times W\times W'$ and
a collection of points $v_1,\ldots,v_m\in D$, $m=nN$, such that
	\begin{enumerate}
	\item the projection $\sigma \colon F\to P^N\times W$ has finite fibres,
	\item for all $f\in Aut_a(D)$ the graph $\Gamma_f$ is a subset
of the closure $\overline{F_{\imath(f)}}$, where the map
$\imath\colon Aut(D)\to P^N$ is given by
$\imath(f) = (v_1,f(v_1),f_{*v_1},\ldots,v_m,f(v_m),f_{*v_m})$.
	\end{enumerate}
\end{Prop}
\end{samepage}

\begin{Rem}\label{any-v}
	Once the set $v_1,\ldots,v_m$ is chosen, we can add to it finitely
many other $v$'s and not change the statement of Proposition~\ref{par1}.
\end{Rem}

	Before we start with the proof we need a technical lemma.

\begin{Lemma}\label{pts}
	Let $A$, $B$, $C\subset A\times B$ be constructible sets and
every fibre $C_a:=\{b\in B\mid (a,b)\in C \}$ be of codimension at least one.
Then there
exists a finite number of points $b_\mu\in B$, $\mu=1,\ldots,s$
such that for every
$a\in A$ there is a point $b_\mu\notin C_a$.
\end{Lemma}

	{\bf Proof.} We first prove the statement for $A$ a locally closed
irreducible subvariety by induction on dimension of $A$.
If $\dim A=0$, the statement is obvious.
Assume it to be proven for $\dim A < d$. By definition of constructible sets,
$C$ is a finite union of locally closed subvarieties
$C_\alpha=U_\alpha \cap F_\alpha$ where $U_\alpha$ are Zariski open and
$F_\alpha$ are closed subvarieties. The subvarieties $F_\alpha$ are not
open, otherwise a fibre $C_a$ would contain an open subset.
So we can choose a point
$$(a_0,b_0)\in (\cap_\alpha U_\alpha) \backslash (\cup_\alpha F_\alpha)$$
The set $A'$ of points $a\in A$, such that $b_0\in C_a$,
is the projection on $A$ of the
intersection $(A\times \{b_0\})\cap C$, which is constructible.
There is an entire neighborhood of $a_0$ in the complement and,
hence, $A'$ has lower dimension than $A$.
Now we use induction for all irreducible components of the
closure $\overline{A'}$. This yields a number of points $b_\mu$. These
points together with $b_0$ satisfy the required condition.

	To prove the statement in case $A$ is constructible we note, that
$A$ is by Definition~\ref{cons} a finite union of locally closed $A_\alpha$'s.
For every $A_\alpha$ with $C_\alpha:=(A_\alpha \times B) \cap C$
the statement of Lemma gives a finite set of points $b_\mu$. The union of
these finite sets for all $\alpha$ satisfies the required property.
\par\hfill {\bf Q. E. D.}

	{\bf Proof of Proposition~\ref{par1}.}
Now we apply Lemma~\ref{pts} to our situation. Let $P'$ be the
constructible subset of all parameters $p\in P$ such that
$E_p\subset V\times W$ is of codimension at least $1$.
Then we set in Lemma~\ref{pts} $A:=P'$, $B:=V\times W$ and
$C:=E \cap (P'\times V\times W)$.
The statement of Lemma yields a number of points
$$(v^{(\mu)},w_\mu) \in D^{n+1},\quad \mu=1,\ldots,N.$$
For every $f\in Aut_a(D)$ and $(v,w)\notin E(f)$ we have by condition~3
in Proposition~\ref{par}, $E(f)\subset E_p$ for $p=(v,f(v),f_{*v})$.
Then for some $\mu=1,\ldots,N$ we have $(v^{(\mu)},w_\mu)\notin E_p$, i.e.
$(v^{(\mu)},w_\mu)\notin E(f)$ and, by condition~2 in
Proposition~\ref{par}, $(p,w,f(w))\in F$. We obtain
$m=Nn$ points $v_1,\ldots,v_m$.

	Now we construct the required family $F$ to be the union of
the sets $F_\mu$ defined by
\begin{equation}\label{fnu}
	F_\mu := \{ (p_1,\ldots,p_N,w,w')\in P^N\times W\times W'
	\mid (p_\mu,w,w')\in F \},
\end{equation}

	Statement~1 in Proposition~\ref{par1} follows from condition~1
in Proposition~\ref{par}. Let $\imath\colon Aut(D) \to P^N$
be the map defined by
\begin{equation}\label{imath}
	\imath(f) := (v_1,f(v_1),f_{*v_1},\ldots,v_m,f(v_m),f_{*v_m}).
\end{equation}
It is in fact a product of Cartan maps $C_v\colon f\mapsto (f(v),f_{*v})$
and is therefore a homeomorphism onto its image.
Statement~2 in Proposition~\ref{par1} follows now from the
above choice of $v_j$'s.
\par\hfill {\bf Q. E. D.}

	\section{Defining conditions for $Aut_a(D)$}

	In Proposition~\ref{par1} we constructed a map
$\imath\colon Aut(D) \to P^N$.
Our goal here is to give semi-algebraic description of the
image $\imath(Aut_a(D))$ and to prove the following Proposition.

\begin{Prop}\label{Ga}
	The image $\imath(Aut_a(D))$ and the set of all graphs
$$ \Gamma := \{(\imath(f),w,f(w)) \mid f\in Aut_a(D) \land w\in D \} $$
are semi-algebraic.
\end{Prop}

\subsection{Reduction to a fixed pattern.}

	In Proposition~\ref{par1} we obtained a constructible family
$F\subset P^N\times W\times W'$. Our goal now is to find
a stratification of $P^N\times W$ such that $F$ has a simplier
form over each stratum. This is done by applying Theorem~\ref{tr} on
local triviality of semi-algebraic morphisms.

	To simplify the notation we shall write $P$ for $P^N$.
We first consider the projection $\sigma\colon F\to P\times W$.
Since we are interested only in points over $P\times D\subset P\times W$,
we write $F\subset P\times D\times W'$ for the intersection with
$P\times D\times W'$. Since $D$ is semi-algebraic, $F$ is semi-algebraic.
The projection $\sigma\colon F\to P\times D$
is a continuous semi-algebraic map (see Definition~\ref{map}) and
we can apply Theorem~\ref{tr} on local triviality.
Theorem~\ref{tr} yields a
finite semi-algebraic stratification $\{Y_1,\ldots,Y_h\}$ of $P\times D$
(see Definition~\ref{str}),
a collection of semi-algebraic typical fibres $\{E_1,\ldots,E_h\}$
and a collection of semi-algebraic structural homeomorphisms
\begin{equation}\label{E}
	\tilde e_i\colon Y_i\times E_i \to \sigma^{-1}(Y_i), \quad i=1,\ldots,h
\end{equation}
(the $\tilde e_i$'s here are the inverses of the $g_i$'s in Theorem~\ref{tr}).
By statement~1 in Proposition~\ref{par1}, every typical fibre $E_i$ is finite.

	The semi-algebraic stratification $\{Y_i\}$ of the product $P\times D$
defines a stratification of every fibre $\{p\}\times D$. This stratification
depends on $p\in P$ and the qualitative picture (e.g. the number of open
strata) can also depend on $p$. However, by changing to a partition of $P$
we reduce this general case to the case of fixed stratification of
$\{p\}\times D$, a fixed {\em pattern}.

	For this we apply Theorem~\ref{tr} again to the projection
$\rho\colon P\times D \to P$ and partition $\{Y_1,\ldots,Y_h\}$ of
$P\times D$. We obtain a finite semi-algebraic stratification
$\{P_1,\ldots,P_r\}$ of $P$, a collection of semi-algebraic typical fibres
$\{G_1,\ldots,G_r\}$, for every $l=1,\ldots,r$ a finite semi-algebraic
partition $\{G_{l1},\ldots,G_{lh} \}$ of $G_l$ and a collection of
semi-algebraic structural homeomorphisms

\begin{equation}\label{G}
	g_l\colon P_l\times G_l \to \rho^{-1}(P_l), \quad l=1,\ldots,r,
\end{equation}
such that
\begin{equation}\label{Gi}
	g_l(P_l\times G_{li}) = (\rho^{-1}(P_l)) \cap Y_i, \quad l=1,\ldots,r,
	\quad i=1,\ldots,h.
\end{equation}

\subsection{The set of all automorphisms.}

	Here we discuss the set of all automorphisms of $D$ the graphs of which
are contained in the closures of fibres of our family $F$.
By condition~1 in
Proposition~\ref{par1}, only finitely many automorphisms can be contained
in the closure of a fixed fibre.
On the other hand, a fixed automorphism can be contained
in closures of a multitude of fibres.

	Without loss of generality we assume, that $\{G_{l1},\ldots,G_{lh} \}$
is a finite semi-algebraic stratification of $G_l$ and $G_{li}$ are connected
(see Proposition~\ref{strat} and Corollary~\ref{con}).

	Now for fixed $p\in P_l$ we wish to determine if the fibre
$F_p\subset D\times W'$ is  related to some $f\in Aut(D)$.
Our procedure for doing this is semi-algebraic:
over the fixed decomposition $D=\sqcup_i G_{il}$
(in fact only over open $G_{il}$'S) we consider the pieces of $F_p$,
determined by the trivialization of it with typical fibres $E_i$.
The condition that certain  of these pieces fit together to form a graph
of an automorphism proves to be semi-algebraic.

	Among the strata $G_{li}$, $i=1,\ldots,h$ we
choose the open one's, which are assumed to be $G_{li}$, $i=1,\ldots,t$,
$t\le h$. By Proposition~\ref{par1},
the projection $\sigma \colon F\to P\times D$ has finite fibres
so the typical fibres $E_i$ are all finite.
Let us fix a $t$-tuple
$e=(e_1,\ldots,e_t)\in E:=E_1\times\cdots\times E_t$.
The number
of possible $t$-tuples is finite. Further, we define the maps $\xi_{e,p}$ over
each open $(Y_i)_p := \{w\in D \mid (p,w)\in Y_i\}\subset D$ by
\begin{equation}\label{xi}
	(p,w,\xi_{e,p}(w)) = \tilde e_i(p,w,e_i), \quad i=1,\ldots,t,
\end{equation}
where $\tilde e_i\colon Y_i\times E_i \to \sigma^{-1}(Y_i)$ are the
trivialization morphisms in~(\ref{E}).

\begin{Prop}\label{pl}
	Let $P_l$ and $e\in E$ be fixed.
The set $P_{e,l}$ of all parameters $p\in P_l$ such that
the map $\xi_{e,p}$ extends to a biholomorphic automorphism from $Aut_a(D)$
is semi-algebraic.
\end{Prop}

	We begin with three lemmas. The first one is a semi-algebraic version of
Lemma~\ref{clos} on constructible sets.

\begin{Lemma}\label{clos-s-a}
	Let $A$, $B$ and $C\subset A\times B$ be semi-algebraic sets. Then the
``fibrewise'' closure of $C$, i.e. the union of closures of the fibres
$C_a :=(\{a\} \times B) \cap C$, $a\in A$ is semi-algebraic.
\end{Lemma}

	{\bf Proof.} Let $X:=A\times B$ and consider the partition
$X_1:=C$, $X_2:=(A\times B)\backslash C$ of $X$. Apply Theorem~\ref{tr}
to the projection of $X$ on $A$. To obtain the
closures of fibres we take the closures of typical fibres $F_{i1}$ in
$F_i$ and their images in $X$ under the structural trivializing
homeomorphisms. The images are semi-algebraic by Theorem~\ref{TS}.
The union of these images for all $i$ is semi-algebraic
and is exactly the ``fibrewise'' closure of $C$.
\nopagebreak\par\hfill {\bf Q. E. D.}

\begin{Lemma}\label{sub}
	Let $A$, $B$ and $C,D\subset A\times B$ be semi-algebraic sets. Then
the set of $a\in A$ such that $C_a \subset D_a$ is semi-algebraic.
\end{Lemma}

	{\bf Proof.} The complement of the required set in $A$ coincides
with the projection on $A$ of the difference $C\backslash D$. The difference
of semi-algebraic set is semi-algebraic, the projection is semi-algebraic
by the Tarski-Seidenberg theorem (Theorem~\ref{TS}).
\par\hfill {\bf Q. E. D.}

\begin{Lemma}\label{1}
	Let $A$, $B$, $C$, $E\subset A\times B$ and $G\subset E\times C$ be
semi-algebraic sets. Then the
set of all $a\in A$, such that for all $b\in E_a$ the fibre $G_{(a,b)}$
consists of exactly one point, is also semi-algebraic.
\end{Lemma}

	{\bf Proof.} We apply the Theorem~\ref{tr} on local trivialization
to the projection of $G$ on $E$. This yields a partition $\{Y_i\}$ of $E$.
The set $E'\subset E$ of one-point fibres $G_{(a,b)}$
is then the finite union of $Y_i$'s
such that the corresponding typical fibres $F_i$ consist of one point.
It follows that $E'$ is semi-algebraic. The required set in $A$ coincides with
the
set of $a\in A$ such that $C_a\subset E'_a$.
The latter set is semi-algebraic by Lemma~\ref{sub}.
\par\hfill {\bf Q. E. D.}

	{\bf Proof of Proposition~\ref{pl}.} We first consider the condition that
$\xi_{e,p}$ extends to a well-defined continuous map on $D$.
This means that for any point
$w\in \overline{(Y_i)_p} \cap \overline{(Y_j)_p}$, $i,j=1,\ldots,t$,
the limits of graphs of $\xi_{e,p}$ over $(Y_i)_p$ and $(Y_j)_p$
coincide over $w$ and consist of one point.

For every stratum
$G_{ls}$, $G_{li}$, $G_{lj}$, $s=1,\ldots,h$, $i,j=1,\ldots,t$, with
$$G_{ls}\subset \overline{(Y_i)_p} \cap \overline{(Y_j)_p}$$
we write these conditions in a form

\begin{equation}\label{cont}
	\left. \matrix{
  i) &  \overline{\Gamma_i(p)} \cap (B(p)\times W') =
	    \overline{\Gamma_j(p)} \cap (B(p)\times W'),       \cr\cr
  ii) & \forall w\in B(p) :
	    \#(\overline{\Gamma_i(p)} \cap (\{w\}\times W')) = 1
	} \right\rbrace,
\end{equation}
where
$$B(p) := g_l(\{p\} \times G_{ls}) \subset \{p\} \times D$$
and
$\Gamma_i(p) := \tilde e_i((Y_i)_p\times \{e_i\})$ is the graph of $\xi_{e,p}$
over $(Y_i)_p$.

	Now, by Lemmas~\ref{clos-s-a} and \ref{sub},
the set $\{p\in P_l\mid i)$ in (\ref{cont}) is satisfied $\}$ is a
semi-algebraic subset of $P_l$. For the condition ii) in (\ref{cont})
we set in Lemma~\ref{1}
$A:=P'$, $B:=D$, $C:= W'$,
$$E:=\{(p,w)\in P'\times D \mid w\in B(p) \}$$
and
$$G:=\{(p,w,w')\in E\times W' \mid
w\in \overline{\Gamma_i(p)} \cap (B(p)\times W')\}.$$
Then, by Lemma~\ref{1}, the set $\{p\in P_l\mid ii)$ in (\ref{cont})
is satisfied $\}$ is a semi-algebraic subset of $P_l$.

Without loss of generality, $i)$ and $ii)$ are satisfied for $p\in P_l$.
Thus, the closures of graphs of $\xi_{e,p}$ over $\{p\} \times D$ yield
well-defined maps $\xi_p\colon D\to  W'$ (We do not know yet, whether or not
these maps are continuous).

	The next condition on $\xi_{e,p}$ is
\begin{equation}\label{D}
	\xi_{e,p}(D) = D,
\end{equation}
which is, by Lemma~\ref{sub}, a semi-algebraic condition.

	Now, if conditions (\ref{cont}) and (\ref{D}) are satisfied,
we can prove that $\xi_{e,p}$ is continuous.
For this let $U_p\subset D$ denote the union of all
$(Y_i)_p$'s, $i=1,\ldots,t$.
This is an open dense subset of $D$ where $\xi_{e,p}$ is continuous.
Fix a point $w_0\in D$. By (\ref{cont}), $\xi_{e,p}(w_0)$ is the only limit
value of $\xi_{e,p}(w)$ for $w\in U_p$. Since $\xi_{e,p}$ is bounded, we have
\begin{equation}\label{lim}
	\xi_{e,p}(w_0) = \lim_{{w\to w_0 \atop w\in U}} \xi_{e,p}(w),
\end{equation}
which means $\xi_{e,p}$ is continuous.

	Thus, we obtained a family of continuous maps $\xi_{e,p}$ from $D$ onto
$D$, which are holomorphic outside some real analytic locally closed
subvariety of codimension~$1$. By the theorem on removable singularities,
$\xi_{e,p}$ is holomorphic on $D$.

	Further, by the theorem of Osgood (see \cite{N}, Theorem~5, Chapter~5)
$\xi_{e,p}$ is biholomorphic if and only if it is injective. This is the
condition on fibres:

\begin{equation}\label{inj}
	\#(\xi_{e,p}^{-1}(y))=1,\, y\in D.
\end{equation}

	The set $\{p\in P_l\mid$ (\ref{inj}) is satisfied $\}$ is,
by Lemma~\ref{1}, semi-algebraic
(we set $A:=P'$, $B:= W'$, $C:=D$, $E:=P'\times D\subset A\times B$ and
$G$ is the family of graphs of $\xi_{e,p}$). This finishes the proof of
Proposition~\ref{pl}.
\par\hfill {\bf Q. E. D.}

\subsection{The image $\imath(Aut_a(D))$ and the set of associated graphs.}

	Here we wish to prove Proposition~\ref{Ga}. Let $P_{e,l}$
be the semi-algebraic subsets from Proposition~\ref{pl}.
We obtain a diagram:
\begin{equation}\label{di}
\def\normalbaselines{\baselineskip20pt\lineskip3pt\lineskiplimit3pt}
\def\mapup#1{\Big\uparrow \rlap{$\scriptstyle#1$}} \matrix{
	     &              & P              \cr
	 & \nearrow & \mapup{\imath} \cr
P_{e,l}  & \to      & Aut_a(D)       \cr
p        & \mapsto  & \xi_{e,p}      \cr },
\end{equation}
where the map from $P_{e,l}$ into $P$ is the usual inclusion.
We define $P'_{e,l}\subset P_{e,l}$
to be the subset of all points $p\in P_{e,l}$,
for which the diagram is commutative. This condition means
$p=(v,v',l)=(v,\xi_{e,p}(v),(\xi_{e,p})_{*v})$ and is therefore semi-algebraic.
Therefore, $P'_{e,l}$ is semi-algebraic. The semi-algebraic property of
$\imath(Aut_a(D))$ is a consequence of the following observation.

\begin{Lemma}
	$$\imath(Aut_a(D))=\bigcup_{{e\in E \atop l=1,\ldots,r}} P'_{e,l}.$$
\end{Lemma}

	{\bf Proof.} Let $p\in \imath(Aut_a(D))$, i.e.
$p=\imath(f)$ for some $f\in Aut_a(D)$. Then, by Proposition~\ref{par1},
$\Gamma_f\subset \overline{F_p}$.
We have $p\in P_l$ for some $l=1,\ldots,r$.
The graph $\Gamma_f$ defines sections in $F$
over every connected open stratum $(Y_i)_p$, $i=1,\ldots,t$. This means
that for some choice $e\in E$ we have $f=\xi_{e,p}$.
Then $\xi_{e,p}\in Aut_a(D)$, which implies $p\in P_{e,l}$. Further, the
equality $f=\xi_{e,p}$ means that diagram~(\ref{di}) is commutative for $p$.
Then $p\in P'_{e,l}$ and the inclusion in one direction is proven.

	Conversely, let $e\in E$ be fixed and $p\in P'_{e,l}$.
Since $p\in P_{e,l}$,
$f:=\xi_{e,p}$ is an automorphism in $Aut_a(D)$.
The commutativity of diagram~(\ref{di}) means
$p=\imath(f)$. This implies $p\in \imath(Aut_a(D))$, which proves the
inclusion in other direction.
\nopagebreak\par\hfill {\bf Q. E. D.}

	{\bf Proof of Proposition~\ref{Ga}.}
The family
$$\Gamma:=\{(\imath(f),w,f(w)) \mid f\in Aut_a(D) \land w\in D \}$$
over $P'_{e,l}$ coincides now with the family of graphs of $\xi_{e,p}$.
The latter is,
by construction, semi-algebraic and Proposition~\ref{Ga} is proven.
\nopagebreak\par\hfill {\bf Q. E. D.}

	\section{Semi-algebraic structures on $Aut_a(D)$}

	Here we finish the proof of Theorem~\ref{main}.
In previous section we considered imbeddings $\imath\colon Aut(D)\to P^N$.
Here we wish to change to Cartan imbeddings $C_v(f):=(f(v),f_{*v})$.
By Proposition~\ref{par1}, $\imath$ is given by
$\imath(f)=(v,f(v),f_{*v})$, where $v=(v_1,\ldots,v_m)\in D^m$.

	The image $C_{v_j}(Aut_a(D))$ is equal to the projection of
$\imath(Aut_a(D))$ on $V'_j\times L_j$ (we use our notations
$v'_j=f(v_j)\in V'_j$, $l_j=f_{*v_j}\in L_j$). By Theorem~\ref{TS} of
Tarski-Seidenberg, $C_{v_j}(Aut_a(D))$ is semi-algebraic. Further, it is
semi-algebraically isomorphic to $\imath(Aut_a(D))$.
By Remark~\ref{any-v}, any $v$ and $v'$ can be
among $v_j$'s. Thus we obtain the following result.

\begin{Prop}\label{stru}
	Let $v\in D$ be any point. The image $C_v(Aut_a(D))$ is semi-algebraic
and this semi-algebraic structure is independent of $v\in D$.
\end{Prop}

	Now we fix some $v\in D$ and denote by $K$ the image
$C_v(Aut_a(D))\subset D\times Gl(n)$. The family
$\Gamma'\subset C_v(Aut_a(D))\times D^2$ of graphs over $C_v(Aut_a(D))$
is a projection of $\Gamma$ and is therefore semi-algebraic.
To simplify our notation we set $P:=D\times Gl(n)$, and $\Gamma:=\Gamma'$.

	Statements~2 and 3 in Theorem~\ref{main} can
now be formulated as follows:

\begin{Lemma}\label{s-a}
	With respect to the group operation of $Aut_a(D)$,
$K$ is a Nash group and the action on $D$ is Nash.
\end{Lemma}

	{\bf Proof.} We consider the graph of the operation
$(x,y)\mapsto xy^{-1}$ in $K^3$. For this, start with the family
$\Gamma\subset K\times D\times D$ and define a new family
$\Gamma_1\subset K^3\times D^3$ by

\begin{equation}\label{G1}
	\Gamma_1 := \{ (x,y,z,w,w',w'')\in K^3\times D^3 \mid
	(y,w',w)\in \Gamma \land (x,w',w'')\in \Gamma \}.
\end{equation}

	The conditions in (\ref{G1}) express the fact that $y^{-1}\in K$
transforms $w$ in $w'$ and $x\in K$ transforms $w'$ in $w''$. The projection
$\Gamma_2$ of $\Gamma_1$ on $K^3\times D^2$ (with coordinates $(x,y,z,w,w'')$)
is, by the Theorem of Tarski-Seidenberg, semi-algebraic. Now the
condition $z=xy^{-1}$ means that the graphs of $z$ and $xy^{-1}$ coincide,
i.e. the fibres $(\Gamma_2)_{(x,y,z)}$ and $(\Gamma_3)_{(x,y,z)}$ coincide,
where $\Gamma_3 := \{(x,y,z,w,w'') \mid (z,w,w'')\in \Gamma \}$ is an
extension of $\Gamma$.

	By Lemma~\ref{sub},     the coincidence of fibres is
a semi-algebraic condition on $(x,y,z)\in K^3$. This proves that the
graph of the correspondence $(x,y)\mapsto xy^{-1}$ is semi-algebraic, which
means that the group operation is semi-algebraic. Since the latter is also
real analytic by Theorem of Cartan, $K$ is an affine Nash group.
Furthermore, the graph $\Gamma$ of the action of $K$ on $D$ is semi-algebraic
and real analytic and therefore Nash.
\par\hfill {\bf Q. E. D.}

	It remains to prove statement~1 in Theorem~\ref{main}
which asserts that $Aut_a(D)$ is a Lie subgroup of $Aut(D)$.

	{\bf Proof of statement~1. } We begin with the semi-algebraic set $K$.
By Proposition~\ref{strat}, it admits a finite semi-algebraic stratification.
Let $x\in K$ be a point in a stratum of maximal dimension. Then there is
a neighborhood $U_x\subset P$ of $x$, such that $K\cap U_x$ is a closed real
analytic submanifold of $U_x$. The preimage
$K':=C_v^{-1}(K\cap U_x)$ is a closed real analytic submanifold in the
neighborhood $U_f:=C_v^{-1}(U_x)$ of $f:=C_v^{-1}(x)$. Since
$Aut_a(D)$ is a subgroup of $Aut(D)$, we see that
$Aut_a(D)\cap (f^{-1} \cdot U_f) = f^{-1} \cdot K'$ is a closed real
analytic submanifold in the neighborhood $f^{-1} \cdot U_f$ of the unit
$id\in Aut(D)$. This implies that $Aut_a(D)$ is a real analytic
subgroup of $Aut(D)$. \hfill {\bf Q. E. D.}


\begin{samepage}
\indent Dmitri Zaitsev

Fakult\"at f\"ur Mathematik

Ruhr-Universit\"at Bochum

44780 Bochum

Germany

\bigskip
e-mail: Dmitri.Zaitsev@ruba.rz.ruhr-uni-bochum.de
\end{samepage}
\end{document}